\documentclass[12pt]{article}
\usepackage{amsmath, amssymb}
\usepackage{enumitem}
\usepackage{authblk}
\usepackage{lineno}
\usepackage[toc,page]{appendix}
\usepackage{amsthm}
\usepackage{float}
\usepackage{changepage}
\usepackage{aliascnt}
\usepackage{theoremref,amsthm}
\usepackage{morefloats}
\usepackage{amsfonts}
\usepackage{amssymb}
\usepackage{graphics}
\usepackage{graphicx}
\usepackage{listings}
\usepackage{xcolor}
\usepackage{float}
\usepackage{verbatim}
\usepackage{geometry}
\usepackage{centernot}
\usepackage[utf8]{inputenc}
\usepackage[english]{babel}
\usepackage{setspace}
\usepackage{thmtools}
\usepackage{xcolor}
\usepackage[bottom]{footmisc}
\usepackage{euscript}
\usepackage{mathtools}
\usepackage{enumitem}
\usepackage [autostyle, english = american]{csquotes}
\MakeOuterQuote{"}
\newcommand{\subscript}[2]{$#1 #2$}
\DeclareMathAlphabet{\itbf}{OML}{cmm}{b}{it}
\providecommand{\keywords}[1]{\textbf{Keywords:} #1}

\newcommand{\bea}{\begin{eqnarray*}}
\newcommand{\eea}{\end{eqnarray*}}
\newcommand{\bean}{\begin{eqnarray}}
\newcommand{\eean}{\end{eqnarray}}
\newcommand{\p}{\partial}
\newcommand{\f}{\frac}

\newcommand{\no}{\nonumber}
\newcommand\ov{\overline}
\newcommand{\ri}{\rightarrow}
\newcommand{\sm}{\setminus}
\newcommand{\aaa}{\mbox{$[$}}
\newcommand{\bbb}{\mbox{$]$}}

\newtheorem{lem}{Lemma}[section]
\newtheorem{prop}{Proposition}[section]

\newtheorem{thm}{Theorem}[section]

\newcommand{\doubleR}{\mathbb{R}}
\newcommand{\RR}{\mathbb{R}}

\renewcommand{\ov}[1]{\overline{#1}}

 \author{ Darko  Volkov 
\thanks{Department of Mathematical Sciences,
Worcester Polytechnic Institute, Worcester, MA 01609. Corresponding author email: darko@wpi.edu. 
} }

\begin{document}
		\title{On the   crack inverse problem for pressure waves in half-space}
		\maketitle

\begin{abstract}
After formulating  the pressure wave
equation  in half-space minus a crack with a zero Neumann condition 
on the top plane,
we introduce a related inverse problem. That inverse problem consists of identifying 
the crack and the unknown forcing term on that crack from overdetermined boundary data 
on a relatively  open set of the top plane.
This  inverse problem is not uniquely solvable unless some additional assumption is made. 
However, we show
that we can differentiate two  cracks $\Gamma_1$ and
$\Gamma_2$ under
the assumption that $\RR^3 \sm \ov{\Gamma_1\cup \Gamma_2}$
is connected. 
If that is not the case we provide counterexamples that demonstrate non-uniqueness, even if   $\Gamma_1$ and
$\Gamma_2$ are smooth and     "almost" flat.
Finally, we show   in the case where  $\RR^3 \sm \ov{\Gamma_1\cup \Gamma_2}$ is not 
necessarily connected that after excluding a discrete set of frequencies, $\Gamma_1$
and $\Gamma_2$ can again be differentiated from overdetermined boundary data.
\end{abstract}


\section{Introduction}
\begin{figure}[H]
    \centering
      \includegraphics[scale=.7]{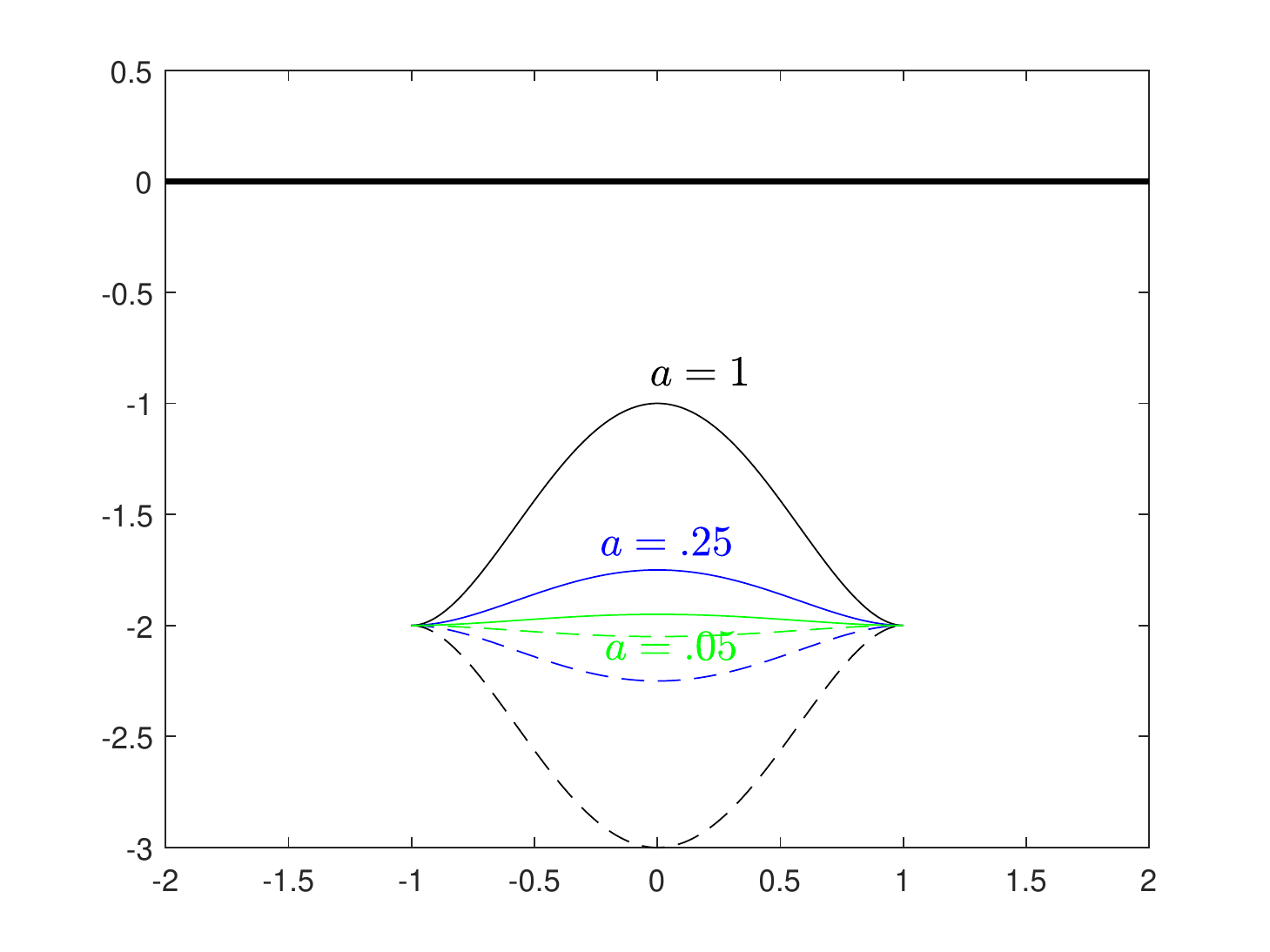}
			       \caption{The shapes $\Gamma_a$ (solid curves) and $\Gamma_a'$ (dashed curves)
						for three values of 
						$a$. We prove that it is impossible to distinguish $\Gamma_a$ from 
						$\Gamma_a'$ at some wavenumbers and some forcing terms on $\Gamma_a$ 
						and $\Gamma_a'$  for the half space crack PDE
						from Cauchy data,
						even if the Cauchy data is given everywhere
						on the top boundary.
}
    \label{three shapes}
\end{figure}

In this paper, we study  an inverse problem consisting of identifying 
a crack in half-space and the unknown forcing term on that crack from overdetermined boundary data 
on a relatively open set of the top plane.
For the forward problem, the governing equations involve the Helmholtz operator,
a zero Neumann condition on the top plane, and continuity of the normal derivative across the crack.
A jump across the crack constitutes the forcing term. Some decay at infinity is enforced  by
requiring that the solution lie in an adequately weighted Sobolev space.\\
Closely related inverse problems have been extensively studied in the steady state case. 
In fact, the steady state case has been investigated for the Laplace operator
\cite{ionescu2006inverse, volkov2021stability}
and for the  linear elasticity operator \cite{aspri2020analysis,
aspri2022dislocations, volkov2017reconstruction,
volkov2017determining, volkov2020stochastic,
triki2019stability, volkov2019stochastic, volkov2022parallel, volkov2022stability, volkov2010eigenvalue}.
	\cite{ionescu2006inverse, volkov2010eigenvalue} cover a related eigenvalue problem derived from 
	stability analysis. 
	In \cite{aspri2020analysis}, 
the direct crack problem for half space elasticity 
was analyzed under weaker regularity conditions (weaker than $H^1$ regular).
 In \cite{aspri2022dislocations},
this direct problem 
was proved to be uniquely solvable in case of
piecewise Lipschitz elasticity coefficients and general elasticity tensors.
Both \cite{aspri2022dislocations}  and \cite{aspri2020analysis}
include a proof of  uniqueness  for the related crack (or fault) inverse problem under appropriate assumptions.	
\cite{triki2019stability}	and \cite{volkov2021stability}  focus  on the Lipschitz stability 
	of the reconstruction of cracks based on the assumption that only planar cracks are admissible.
	The seismic model  introduced in 
	\cite{volkov2017reconstruction}
	was used in
	\cite{volkov2017determining}
	to address a real life problem in geophysics consisting of identifying a fault 
	using GPS measurements of surface displacements.
	\cite{volkov2020stochastic, volkov2019stochastic}
	feature statistical numerical methods for the reconstruction of cracks 
	based on a Bayesian approach.
	In \cite{volkov2022parallel}, a related parallel accept/reject sampling 
	algorithm was derived.
	The numerical method in
 \cite{volkov2022stability} is entirely different, it is based on deep learning.\\
Here, we analyze  a direct and an inverse problem that generalize the Laplace based case
to the wave case. The waves considered here are time harmonic pressure waves and can
 be modeled by an inhomogeneous Helmholtz equation.
 We can actually model heterogeneous media by 
assuming that the 
wavenumber $k^2$ is an $L^\infty$ function, as long as it is non-negative, bounded away from zero, and
equal to a constant $k_0$ outside a bounded set.
In section 
	\ref{direct and inverse crack problems}, we prove that the direct pressure wave
	problem in half space minus a crack is uniquely 
	solvable and well posed in the 
	space of functions 
\bea\{ v \in H^1_{loc}(\RR^{3-}\setminus \ov{\Gamma}): \f{v}{\sqrt{1+r^2}}, 
\f{\nabla v}{\sqrt{1+r^2}}, \f{\p v}{\p r} - i k_0 v \in L^2(\RR^{3-}\setminus \ov{\Gamma})
\},
\eea
where $\Gamma$ is the crack and  $r$ is the distance to the origin.
The proof of uniqueness for the inverse problem relies on unique continuation for elliptic equations
from Cauchy data. 
The earliest such continuation results relied on properties of analytic functions.
Later, Nirenberg proved that 
for second order PDEs whose leading term is the Laplacian,
it suffices to assume that solutions are $C^1$ with piecewise continuous second derivatives
 for the unique continuation property to hold \cite{nirenberg1957uniqueness}.
This result was further improved by Aronszajn et al. 
where in \cite{aronszajn1962unique}
it was extended to such PDEs with only Lipschitz coefficients.
Unfortunately, demanding Lipschitz continuity is impractical in applications since it does not 
even cover the piecewise constant case.
More recently, Barcelo et al. \cite{barcelo1988weighted}
proved a stronger result. 
In particular, their unique continuation result implies that 
a solution to the pressure wave  equation $(\Delta + k^2) u =0$ in an open
set of $\RR^3$ satisfies the unique continuation property
if $k^2$ is in 
$L^\infty_{loc}(\RR^3)$.
This unique continuation property will help us show in section \ref{crack inverse problem uniq}
 a uniqueness result for the  pressure wave inverse problem in the half-space 
$\{x: x_3 <0\}$
 minus a crack: 
	if $\Gamma_i$, $i=1,2$
	are two cracks where the forcing terms $g_i$, $i=1,2$ defining pressure discontinuity 
	across $\Gamma_i$
	have full support in $\Gamma_i$ leading to the solution $u_i$ of the forward problem, if
	$\RR^3 \sm \ov{\Gamma_1\cup \Gamma_2}$
is connected, and the Cauchy data for $u_1$ and $u_2$ 
are the same on a  relatively open set of the top boundary
 $\{x: x_3 =0\}$ then $\Gamma_1= \Gamma_2$ and $g_1=g_2$.\\
In section 
\ref{two counter examples}, 
	we show counterexamples where uniqueness for the crack inverse problem
	fails if $\RR^3 \sm \ov{\Gamma_1\cup \Gamma_2}$
is not connected.
In a first class of counterexamples, $\ov{\Gamma_1\cup \Gamma_2}$
is a sphere and we use the first Neumann eigenvalue for the Laplace operator inside an open ball
that is odd about the equator. Such a function is necessarily zero on the equator 
and thus the values on the top half sphere
can be extended by zero to the lower half sphere without losing its $H^{\f12}$ 
character.
One might argue that this first counterexample is unsatisfactory in the sense that it involves a geometry
that is quite "round": if $n$ is the normal vector on $\Gamma_1$ pointing up
the range of $n \cdot e_3$ is $(0,1]$. 
For that reason we provide a second counterexample where that range can be made arbitrarily
narrow, and such that $\Gamma_1$ can be continued into a plane
in a $C^1$ regular fashion. 
Constructing this family of counterexamples requires using arguments borrowed
from the analysis of elliptic PDEs on domains
with cusps. To make our argument easier to follow, this example is 
first constructed in 
dimension 2, then generalized to the three dimensional case
using cylindrical coordinates. 
Let $\Gamma_a$ be the open curve $x_2=a(x_1-1)^2(x_1+1)^2 - 2$, $x_2 
\in (-1,1)$, where $a>0$ is a flattening parameter. 
In figure  \ref{three shapes} we sketched $\Gamma_a$ for $a=1, .25, .05$.
We also sketched in the same figure $\Gamma_a'$, obtained from $\Gamma_a$
by symmetry about the line 
$x_2=-2$. We show that for some values of $k$ and some choice of
forcing terms $g$ on $\Gamma_a$  and $g'$ and $\Gamma_a'$, 
the half space crack PDE leads to the same Cauchy data everywhere on the top boundary 
$x_2=0$. Next, this geometry is rotated in three dimensional space,
and we show how to construct a counterexample to uniqueness using
rotationally invariant forcing terms $g_1, g_2$. 
The values of $g_1, g_2$ on a cross-section are not those from the two dimensional case: slight adjustments have to be made as 
the volume element in cylindrical coordinates is $rdrd \theta d x_3$
and the Laplace operator for a function which is independent of 
$\theta$ is $\f{1}{r} \f{\p }{\p r} r  \f{\p }{\p r} + \p_{x_3}^2$.\\
In section \ref{unique discrete}, our last result states that the pressure wave crack inverse problem 
is uniquely solvable within the class 
${\cal P}$ of Lipschitz open surfaces that are finite unions of polygons, except possibly 
for a discrete set of frequencies. 
A change of frequency amounts to changing the wavenumber $k^2$ to $t^2k^2$,
for some $t>0$.
If $\Gamma_1$ and $\Gamma_2$ 
are in ${\cal P}$, and if $t$ is not in some discrete set, we show that  
if the values of the corresponding solutions $u_1, u_2$  to the 
pressure wave  crack inverse problem in half space  for the wavenumber $t^2 k^2$,
are  equal on a relatively open set of the top boundary 
then $\ov{\Gamma_1} = \ov{\Gamma_2}$, as long as the jump of $u_i$ has full support 
across $\Gamma_i$.
	
\section{The direct and inverse crack problems 
for pressure waves in half space}\label{direct and inverse crack problems}
\subsection{Problem formulation and physical interpretation} 
Let
$\RR^{3-}$ be the open half space  $\{x=(x_1,x_2,x_3):  x_3<0 \}$. 
Let $\Gamma$ be a Lipschitz  open surface in $\RR^{3-}$.
Assume that $\Gamma$ is strictly included in $\RR^{3-}$ so that  the distance 
	from $\Gamma$ to the plane $\{ x_3 =0 \}$ is positive and let
	$D$ be a domain in $\RR^{3-}$
	with Lipschitz boundary such that $\Gamma \subset \p D$.
The trace
	theorem (which is also valid in Lispchitz domains, \cite{ding1996proof, gagliardo1957caratterizzazioni}),
	allows us to define an inner and outer trace in $H^{\f12}(\p D)$ 
	of functions defined in $\RR^{3-} \setminus \p D$ with local $H^1$ regularity.
Let $\tilde{H}^{\f12} (\Gamma)$ be the space of functions
	in $H^{\f12}(\p D)$ supported in $\ov{\Gamma}$.
	Let $k$ be in $L^\infty(\doubleR^{3-})$  such that,
\begin{enumerate}[label=(\subscript{H}{{\arabic*}})]
\item $k$ is real-valued,
\item there is a positive constant $k_{min}$ such that $k \geq k_{min}$
almost everywhere in $\doubleR^{3-}$,
\item there is an $R_0>0$ and a $k_0>0$ such that if $|x|\geq R_0$, and $x \in 
\doubleR^{3-}\sm \ov{\Gamma}$,
$k(x) = k_0$.
\end{enumerate}
We define the direct pressure wave crack problem in half-space to be the boundary value problem, 
	\bean
        (\Delta + k^2 )u=0\text{ in }\doubleR^{3-}\sm \ov{\Gamma},  \label{BVP1}     \\
        \p_{x_3} u=0\text{ on  the surface } x_3=0,  \label{BVP2}     \\
     \aaa \frac{\p u}{\p n} \bbb = 0 \mbox{ across }\Gamma,   \label{BVP3}\\
      \aaa u \bbb
			=g \mbox{ across }\Gamma,  \label{BVP4}\\
    u \in {\cal V} \label{Decay1},
\eean
where  $n$ is a unit normal vector  on $\Gamma$, 
$\aaa v \bbb $ denotes the jump of a function $v$ across $\Gamma$ in the normal
direction,
$g$ is in $\tilde{H}^{\f12} (\Gamma)$,
and ${\cal V}$ is a functional space ensuring that this PDE is well posed and that 
the solution $u$ depends continuously on $g$.
Following \cite{nedelec2001acoustic}, section 2.6, we introduce
\bean
{\cal V} := \{ v \in H^1_{loc}(\RR^{3-}\setminus \ov{\Gamma}): \f{v}{\sqrt{1+r^2}}, 
\f{\nabla v}{\sqrt{1+r^2}}, \f{\p v}{\p r} - i k_0 v \in L^2(\RR^{3-}\setminus \ov{\Gamma})
\},
\eean
where as previously $r=|x|$ for $x $ in $\RR^3$. 
\begin{prop} \label{direct  prob}
Problem (\ref{BVP1}-\ref{Decay1}) is uniquely solvable.
The solution 
$u$ in ${\cal V}$  depends continuously on the forcing term
$g$ in $\tilde{H}^{\f12} (\Gamma)$.
\end{prop}
\textbf{Proof:}
We first show uniqueness. Assume that $g=0$.
It then follows from (\ref{BVP1}-\ref{Decay1}) that $(\Delta + k^2) u = 0$, weakly in
$\RR^{3-}$. Next we extend $u$ to $\RR^3$ by setting
$u(x_1,x_2, x_3) = u(x_1,x_2, - x_3)$ if $x_3 >0$ and $k$ is extended in a similar fashion.
 Since $u \in {\cal V}$, using 
\eqref{BVP2}  we have that $(\Delta + k^2) u = 0$, weakly in
$\RR^{3}$.
Let $S_R$ be the sphere centered at the origin with radius $R$. 
 Applying Green's theorem shows that
$\mbox{Im} \int_{S_R} u \f{\p \ov{u}}{\p r} = 0$.
Next, since $u \in {\cal V}$, 
there is a sequence $R_n \ri \infty$ such that,
\bea
\lim_{n \ri \infty}\int_{S_{R_n}} |\f{\p u}{\p r} - i k_0 u|^2 = 0,
\eea 
so altogether we have that,
\bea
\lim_{n \ri \infty}\int_{S_{R_n}} |\f{\p u}{\p r}|^2 + | u|^2 = 0.
\eea 
Due to Rellich's lemma for far field patterns, it follows that $u (x) =0$, if $|x| >R_0$. 
Since the only regularity assumption on  $k$ is that it is in $L^\infty$ 
there is no elementary argument for showing that $u(x)=0$ if $|x| \leq R_0$.
However, we can use results from the unique continuation literature.
As in our case $k^2$ is in $L^\infty$, 
$u$ is in $H^2_{loc}$ so
the corollary of  theorem 1 in \cite{barcelo1988weighted} 
can be used 
to claim that $u$ is zero throughout $\RR^3$ since 
$2 >\f{6n-4}{3n+2}$, $n=3$.\\
Next, we prove existence of a solution to
(\ref{BVP1}-\ref{Decay1}).
Fix $R' >R_0$. Let $B_{R'}$ the open ball centered at the origin of $\RR^3$
with radius $R'$
and define the bilinear functional,
\bea
{\cal B}(v,w) = \int_{B_{R'}} \nabla v \cdot \nabla w -  k^2  v w  - \int_{S_{R'}} T_{R', k_0} v w,
\eea
for $ v, w \in H^1(B_{R'}) $ and where $T_{R',k_0}$ is the Dirichlet to Neumann map for radiating solutions to the Helmholtz equation in the exterior of $B_{R'}$
with wavenumber $k_0$. $T_{R',k_0}$ is known to be a  continuous
mapping from $H^{\f12 } (S_{R'})$ to $H^{-\f12 } (S_{R'})$, 
while $-T_{R',0}$ is strictly coercive, and $T_{R',k_0} - T_{R',0}$ is compact
from $H^{\f12 } (S_{R'})$ to $H^{-\f12 } (S_{R'})$, see 
\cite{colton1998inverse}, section 5.3, or \cite{nedelec2001acoustic}, section 2.6.5.
According to the uniqueness property covered above, we have that if
$v \in H^1(B_{R'})$ and
${\cal B}(v,w) =0 $ for all $w \in H^1(B_{R'})$, then $v=0$. \\
Let $u_g$ be in $H^1(B_{R'} \setminus \ov{\Gamma})$ 
such that $\aaa u_g \bbb = g$ across $\Gamma$ and $u_g$ is zero in a neighborhood of
$S_{R'}$. For example, we can set
\bea
u_g(x) = \phi(x) \f{1}{4 \pi}\int_\Gamma \nabla_y (\f{1}{|x-y|}) \cdot n(y) g(y) d \sigma(y),
\eea
where $\phi \in C^\infty_c (B_{R'})$ is constantly equal to 1 in a neighborhood of 
$\ov{\Gamma}$.
It is known that the $H^1$ norm of $u_g$ is bounded by a constant times the 
$H^{\f12}$ norm of $g$, see theorem 1 in \cite{costabel1988boundary}.
Consider the variational problem,
\bean \label{var pb}
\mbox{find } v  \in H^1(B_{R'}) \mbox{ such that } \forall w \in H^1(B_{R'}), \no \\
{\cal B}(v,w) = -{\cal B}(u_g,w).
\eean
We already know that this problem has at most one solution. 
Existence follows by arguing that this problem is in the form strictly coercive plus compact,
which is the case thanks to the properties of the operator $T_{R',k_0}$ recalled above.
Next, integrating by parts will show that $(\Delta + k^2) (v+u_g) =0$, in 
$B_{R'} \setminus \ov{\Gamma}$  and that  $ [ \f{\p }{\p n } (v + u_g)] =0 $ across $\Gamma$.
$ [  (v +u_g)] =[u_g] = g $ across $\Gamma$ is clear by construction. 
Thanks to the operator $T_{R', k_0}$,
$v$ can be extended to a function in $H^1_{loc} (\RR^3 \sm \ov{\Gamma} )$
such that 
$\f{v}{\sqrt{1+r^2}}, 
\f{\nabla v}{\sqrt{1+r^2}}, \f{\p v}{\p r} - i k_0 v\in L^2(\RR^{3}\setminus \ov{\Gamma})$
and $(\Delta + k^2) v =0 $ in $\RR^3\sm  \ov{ B_{R'}}$ as well as  in a neighborhood of $S_{R'}$.
$u_g$, which is zero in a neighborhood of $S_{R'}$ is just extended by zero.
Finally, we set for $x_3 <0$
\bea
u(x_1, x_2, x_3) =
(v+ u_g) (x_1, x_2, x_3)+ 
(v+ u_g)(x_1, x_2, -x_3) 
 \eea
to find a solution to (\ref{BVP1}-\ref{Decay1}). $\square$ 

\subsection{The   crack inverse problem for pressure waves in half space: formulation and uniqueness of solutions}\label{crack inverse problem uniq}
We now  prove a theorem stating that the crack inverse problem related to
problem (\ref{BVP1}-\ref{Decay1}) has at most one solution.
The data for the inverse problem is Cauchy data over a portion
of the top plane $\{ x_3 = 0 \}$.  The forcing term $g$ and the crack $\Gamma$ are both unknown 
in the inverse problem.  

\begin{thm}
\label{InverseProblemResult}
For $i=1,2$,
    let $\Gamma_i$   be a Lipschitz open surface
			such that its closure is
			in
		 $\doubleR^{3-}$, let 
 $u^i$ be the unique solution to  (\ref{BVP1}-\ref{Decay1}) with $\Gamma_i$ in place of $\Gamma$ and the jumps $g^i$ in  $\tilde{H}^{1/2}(\Gamma_i)$ in place of $g$.  
	Let $V$ be a non-empty relatively open subset of the top plane $\{x: x_3=0\}$. 
	Assume that  $\RR^{3-} \setminus \ov{\Gamma_1 \cup \Gamma_2}$ is
	connected and that $g^i$ has full support in  $\ov{\Gamma_i}$, $i=1,2$. 
	If $u^1=u^2$ on $V$, then 
	$\ov{\Gamma_1}=\ov{\Gamma_2}$ and $g^1=g^2$ almost everywhere.	
\end{thm}

\textbf{Proof}:
Let $U=\RR^{3-} \setminus \ov{\Gamma_1 \cup \Gamma_2}$ and set $u= u^1 - u^2$ in $U$.
Since $(\Delta + k^2) u =0$ in $U$ and $ u = \p_{x_3} u = 0$ on $V$, $u$ can be extended by zero
 to an open set $U'$ of $\RR^3 $ such that  $U \subset U'$, $U'$ is connected, 
$U' \cap \{x: x_3>0\}$ is non-empty, $u$ is in $H^2_{loc} (U')$ and satisfies $(\Delta + k^2) u  =0$ in 
$U'$, where we can set $k=0$ in  $U' \cap \{x: x_3>0\}$. 
As $u$ is in $H^2_{loc} (U')$ and $k^2$ is in 
$L^\infty (U')$,
  the unique continuation property (corollary of  theorem 1 in \cite{barcelo1988weighted}) can be applied,
 and $u$ is zero in $U'$. \\
Arguing by contradiction, suppose that there is an $x$ in $\Gamma_1$ such that 
$x \notin \ov{\Gamma_2}$. Then there is an open ball $B(x,r)$ centered at $x$
with radius $r>0$ such that $B(x,r) \cap \ov{\Gamma_2} = \emptyset$. Since $\aaa u \bbb=
\aaa u^2 \bbb =0$ across $B(x,r) \cap \ov{\Gamma_1} $, it follows that
$\aaa u^1 \bbb  =0 $ across $B(x,r) \cap \ov{\Gamma_1} $: this contradicts
that $g^1$ has full support in  $\ov{\Gamma_1}$. We conclude that 
$\Gamma_1 \subset \ov{\Gamma_2}$. Reversing the roles of $\Gamma_1$ and 
$\Gamma_2$ we then find that $\ov{\Gamma_1} =\ov{\Gamma_2}$. Using one more time
that $u$ is zero in $U$,
since $[u]=0$ across $\Gamma_1= \Gamma_2$, it follows that
$g_1- g_2 =0$ almost everywhere in $\Gamma_1$.
$\square$.

\section{Examples where $ \RR^{3-} \setminus \ov{\Gamma_1 \cup \Gamma_2}$ is not connected and
uniqueness for the crack inverse problem for pressure waves fails}\label{two counter examples}

\subsection{A counter example involving half-spheres}\label{counter-example}
We start from an eigenvalue for the Neumann problem in the ball $B(0,1)$, $k_1>0$,  and  $\psi$
an associated eigenfunction 
\bea
(\Delta + k_1^2) \psi =0, \mbox{ in } B(0,1), 
\f{\p \psi}{ \p n} = 0, \mbox{ on } B(0,1). 
\eea
We may choose $\psi$ to be odd in $x_3$ so that $\psi(x_1, x_2, 0) =0$.
More specifically, 
denoting $j_1(r)= \f{\sin r}{r^2} - \f{\cos r}{r}$ the spherical harmonic of order 1,
it is known that $s(r,\theta,\phi) = j_1(r) \cos \theta$
(where in the spherical coordinates $(r, \theta, \phi)$,
$\theta$ is the co-latitude ) satisfies  $ (\Delta +1 ) s = 0$.
Now let $k_1$ be the first positive zero of $j_1'$.
Let $\psi(r, \theta, \phi) = j_1(k_1 r) \cos \theta$, and let $S(0,1)$ 
be the unit sphere centered at the origin. Define the half spheres
$S^+(0,1)=\{x \in S(0,1): x_3>0\}$, $S^-(0,1)=\{x \in S(0,1): x_3<0\}$. \\
Now we set
\bea
\Gamma_1 = \{ x \in \RR^{3-}: x + (0,0,2)  \in S^+(0,1) \}, \\
g_1(x_1,x_2,x_3) = \psi(x_1,x_2,x_3+2), (x_1,x_2,x_3) \in \Gamma_1,
\eea
and likewise,
\bea
\Gamma_2 = \{ x \in \RR^{3-}: x + (0,0,2)  \in S^-(0,1) \}, \\
g_2(x_1,x_2,x_3) = \psi(x_1,x_2,x_3+2), (x_1,x_2,x_3) \in \Gamma_2.
\eea
We decide that on $\Gamma_1 $ and on $\Gamma_2$, the normal vector will have a positive 
third coordinate (pointing "up").
Note that by construction $g_1 \in \tilde{H}^\f12 (\Gamma_1)$ 
and $g_2 \in \tilde{H}^\f12 (\Gamma_2)$ since  $\psi(x_1, x_2, 0) =0$.
Let 
 $u^i$ be the unique solution to  (\ref{BVP1}-\ref{Decay1}) with $\Gamma_i$ in
 place of $\Gamma$ and the jumps $g^i$ in  $\tilde{H}^{1/2}(\Gamma_i)$ in place of $g$.
Let $u= u^1 - u^2$.\\
We now show that $u$ is zero outside $B(0,1) - (0,0,2)$.
Set
\bea
v(x_1,x_2,x_3) = \left\{
\begin{array}{l}
\psi(x_1,x_2,x_3+2), \mbox{ if } x \in B(0,1) - (0,0,2),\\
0, \mbox{ if } x \in \RR^{3-}\sm \ov{B(0,1) - (0,0,2)}.
\end{array}
\right.
\eea
By construction, $[ \f{ \p v}{\p n}]=0$, and $ [v] = \psi = [u]$.
Since problem (\ref{BVP1}-\ref{Decay1})  is uniquely solvable,
$ u =v$ in $\RR^{3-}$. We conclude that $u^1=u^2$ everywhere on the top plane
$\{x: x_3=0\}$. Note that $g_i$ has full support on $\Gamma_i$, $i=1,2$.\\
By a simple rescaling argument, for any positive constant $k^2$, we can likewise construct two 
half-spheres
$\Gamma_1, \Gamma_2$ and forcing terms 
$g_1, g_2$ with full support in  $\Gamma_1, \Gamma_2$ 
such that corresponding $u^1, u^2$ 
solving (\ref{BVP1}-\ref{Decay1})  satisfy  $u^1 = u^2$ everywhere 
on the top plane $\{x: x_3 = 0 \}$.

\subsection{A counterexample  in case of cracks
that are nearly flat, the two-dimensional case}\label{counter-example2}
\begin{figure}[htbp]
    \centering
      \includegraphics[scale=.7]{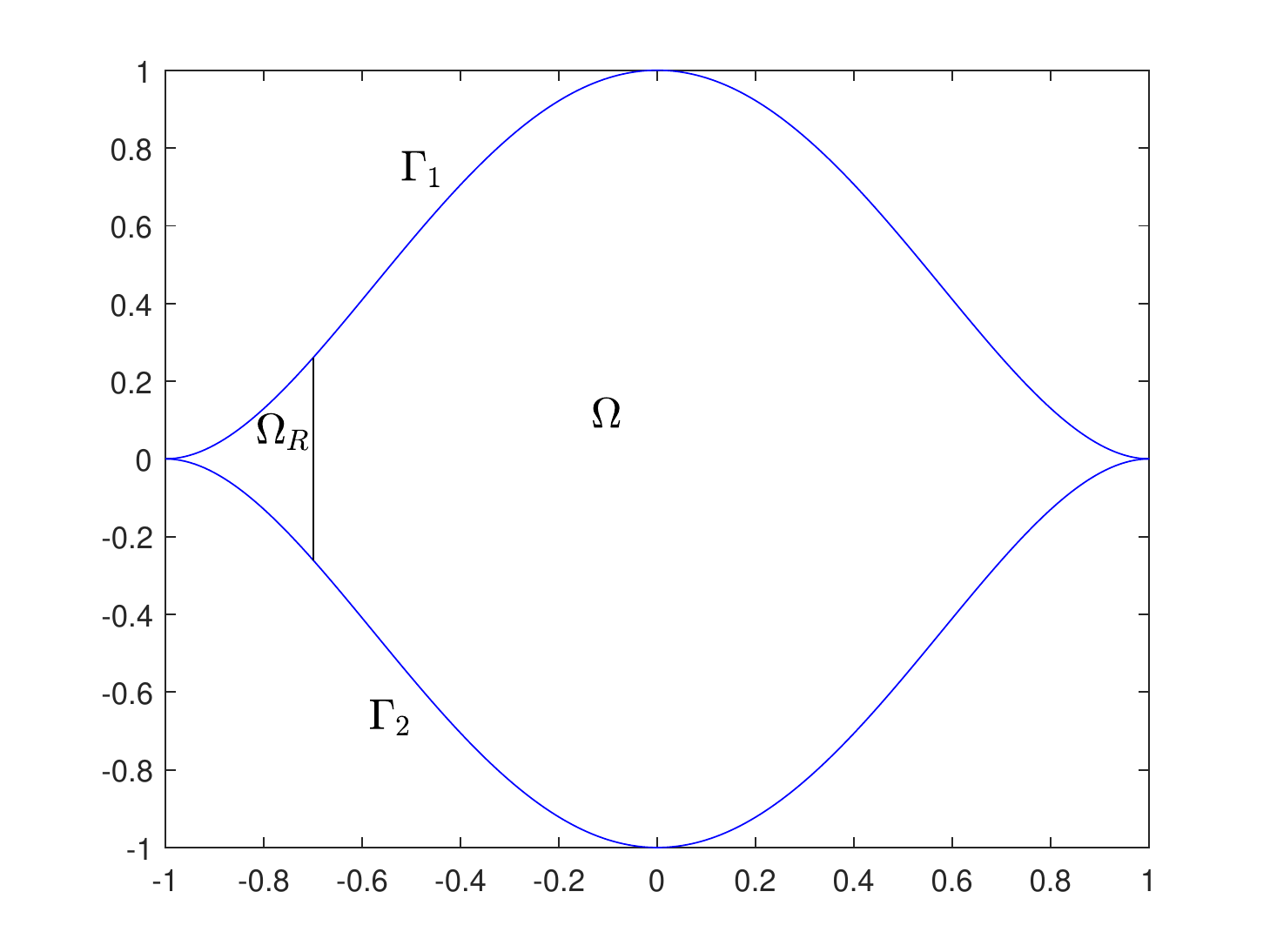}
			       \caption{The domain $\Omega$ bounded by $\Gamma_1$ and
						$\Gamma_2$ and the subdomain $\Omega_R$.
}
    \label{omegaRdomain}
\end{figure}

In the previous counterexample
the range of $n \cdot e_3$, where $n$ is the normal vector to $\Gamma_1$ and $e_3=(0,0,1)$,   was wide, namely 
this range was $(0,1]$.
By contrast, in this section we construct a counterexample where this range can 
be made arbitrarily small. 
The price to pay is that 
constructing this counterexample requires a thorough 
regularity analysis on domains with cusps. We will construct this counterexample in 
dimension 2 to make 
the argument clearer and the notations  easier to follow. 
Let $f:[-1,1] \ri \RR$, $f(t)=(t-1)^2 (t+1)^2$, and consider the domain with two cusps
\bean \label{omega def}
\Omega :=\{ (x_1,x_2):  -1<x_1<1, -f(x_1) < x_2 < f(x_1) \}.
 \eean
$\Omega$ is sketched in figure \ref{omegaRdomain}.

\begin{lem}\label{comp embed}
$H^1(\Omega) $ is compactly embedded in $L^2(\Omega)$.
\end{lem}
\textbf{Proof:} $\Omega $ presents two power-type cusps. This lemma 
holds thanks to the theory of Sobolev spaces on domains with cusps. 
More precisely, we refer to Maz'ya and Poborchi's textbook, p 430, \cite{maz1998differentiable}.
$\square$\\
Thanks to lemma \ref{comp embed}, classical arguments can be used to show that 
the eigenvalues for the Neumann problem for the Laplace operator in $\Omega$ 
form an increasing sequence  diverging to infinity. We now claim that we can find an
eigenfunction which is odd in $x_2$:
\begin{lem} \label{odd u}
There exists a  function $\varphi$ in $H^1(\Omega)$ and $\mu>0$ such that 
$\int_\Omega \varphi^2=1$,
$(\Delta+\mu^2) \varphi =0$, $\f{\p \varphi}{\p n} =0$, almost everywhere
on $\p \Omega$ and $\varphi(x_1, -x_2) = - \varphi(x_1, x_2)$, for $(x_1,x_2)$
  in $\Omega$.
\end{lem}
\textbf{Proof:}
Let
\bea
H^{1,\pm}(\Omega):= \{v \in H^1(\Omega): v(x_1,-x_2) = \pm v(x_1,x_2), (x_1,x_2) 
\in \Omega \}.
\eea
$H^{1,+}(\Omega)$ and $H^{1,-}(\Omega)$ are orthogonal complements of each other
 with regard to the inner 
product  in $H^1(\Omega)$.
Let
\bea
\mu^2= \inf_{v \in H^{1,-}(\Omega), \int_\Omega v^2 =1} \int_\Omega |\nabla v|^2.
\eea
Standard arguments show that a minimizing sequence for this inf converges to some $\varphi$ in
$H^{1,-}(\Omega)$ such that for all $\phi$ in $\in H^{1,-}(\Omega)$
$\int_\Omega \nabla \varphi \cdot \nabla \phi - \mu^2 \varphi \phi  =0$. 
Now since $\varphi \in H^{1,-}(\Omega)$, for all $\phi$ in $\in H^{1,+}(\Omega)$
$\int_\Omega \nabla \varphi \cdot \nabla \phi - \mu^2 \varphi \phi  =0$. 
We conclude that $(\Delta + \mu^2) \varphi =0$  and
$\f{\p \varphi}{\p n} =0$, almost everywhere
on $\p \Omega$. $\square$\\
We would now like to use $\varphi$ to construct a relevant counterexample
just in the way that $\psi$ was used in section \ref{counter-example}. 
However, one important point remains. Since $\Omega$ presents cusp singularities,
not all functions in $H^1(\Omega)$ can be extended to functions in $H^1(\RR^2)$.
The following analysis will show that such an extension is possible for $\varphi$.
The argument borrows from the theory of elliptic PDEs on domains with cusps
and takes advantage of the fact that $\varphi(x_1, 0) =0$, $-1< x_1 <1$, in the sense of traces.
\begin{lem} \label{odd u prop}
Let $\varphi$ be as in lemma \ref{odd u}.
Let $R \in (0,1)$ and 
$$\Omega_R :=\{  (x_1,x_2) \in \Omega: -1 < x_1 < -1 +R\}, $$
(see figure
 \ref{omegaRdomain}).
Let $\sigma$ be the weight $\sigma(x_1,x_2) = \sqrt{(x_1+1)^2 + x^2 }$.
Then
\bean
\int_{\Omega_R} \varphi^2 = O(R^4), \label{OR4}
\eean
and $\sigma^{- \alpha} \varphi$ is in $L^2(\Omega)$ for any $\alpha <2$.
\end{lem} 
\textbf{Proof:}
Let $\eta \in C^1(\Omega) \cap H^1(\Omega)$ be such that  $\eta(x_1,0) =0$.
For $0 < x_2 < f(x_1)$,
\bea
|\eta(x_1,x_2)| & =& |\int_0^{x_2} \p_{x_2} \eta(x_1, t) dt|\\
& \leq & |x_2|^{\f12}  (\int_0^{x_2} |\p_{x_2} \eta(x_1, t)|^2 dt)^{\f12},
\eea
thus we have for $R \in (0,1)$, 
\bea
\int_{-1}^{-1+R} \int_0^{f(x_1)}
|\eta(x_1,x_2)|^2 d x_2 d x_1 & \leq & 
\int_{-1}^{-1+R}\int_0^{f(x_1)}
(|x_2| \int_0^{x_2} |\p_{x_2} \eta(x_1, t)|^2 dt)
d x_2 d x_1 \\
&\leq &  \int_{-1}^{-1+R} \int_0^{f(x_1)}
(|x_2| \int_0^{f(x_1)} |\p_{x_2} \eta(x_1, t)|^2 dt)
d x_2 d x_1 \\
&\leq &  \int_{-1}^{-1+R} \f{f(x_1)^2}{2} \int_0^{f(x_1)} |\p_{x_2} \eta(x_1, t)|^2 dt d x_1 \\
&\leq & 8 R^4  \int_{-1}^{-1+R} \int_0^{f(x_1)} |\p_{x_2} \eta(x_1, t)|^2 dt d x_1 ,
\eea
that is, $\int_{\Omega_R} |\eta|^2 \leq  16  R^4 \int_{\Omega_R} |\p_{x_2} \eta|^2 $.
Estimate \eqref{OR4} follows by density since
$\varphi(x_1,0)=0, -1<x_1<1$, in the sense of traces. 
Next,
\bea 
&&\int_{\Omega_1} |\varphi(x_1, x_2)|^2 ((x_1+1)^2 + x^2)^{- \alpha} dx_1 dx_2\\
&=& \sum_{n=0}^\infty \int_{-1+\f{1}{2^{n+1}}}^{-1+\f{1}{2^{n}}}
\int_{-f(x_1)}^{f(x_1)} |\varphi(x_1, x_2)|^2 ((x_1+1)^2 + x^2)^{- \alpha} dx_2 dx_1\\
&\leq & C \sum_{n=0}^\infty \f{2^{(2n+2) \alpha}}{2^{4n}},
\eea
is finite if $\alpha <2$, thus $\sigma^{- \alpha} \varphi$ is in $L^2(\Omega)$, if $\alpha <2$.
$\square$ \\
\begin{lem} \label{grad u prop}
Let $\varphi$ be as in lemma \ref{odd u}.
$\sigma^{- \beta} \nabla \varphi$ is in $L^2(\Omega)$, if $\beta <1$.
\end{lem}
\textbf{Proof:}
We know that 
$\int_\Omega \nabla \varphi \cdot \nabla \phi - \mu^2 \varphi \phi  =0$, for all $\phi$
in $H^1(\Omega)$.
Let $\beta <1$, $0< \epsilon <1$, and
\bea
\phi_\epsilon = \left\{ \begin{array}{l}
\varphi \sigma^{-2 \beta}, \mbox{ if } \sigma> \epsilon \\
\varphi \epsilon^{-2 \beta}, \mbox{ if } \sigma \leq \epsilon
\end{array}
\right.
\eea
$\phi_\epsilon$ is in $H^1(\Omega)$.
From $\int_\Omega \nabla \varphi \cdot \nabla \phi_\epsilon = \mu^2 \varphi \phi_\epsilon $,
\bea
&& \int_\Omega (|\nabla \varphi|^2 \sigma^{-2 \beta} - 2 \beta
\varphi \sigma^{-2 \beta -1} \nabla \varphi \cdot \nabla\sigma) 1_{\sigma > \epsilon}
+ \int_\Omega |\nabla \varphi|^2 \epsilon^{-2 \beta} 1_{\sigma \leq \epsilon} \\
 &= & \mu^2 (
\int_\Omega \varphi^2 \sigma^{-2 \beta}  1_{\sigma > \epsilon}
+ \int_\Omega  \varphi^2 \epsilon^{-2 \beta} 1_{\sigma \leq \epsilon} ),
\eea
we infer,
\bea
\int_\Omega |\nabla \varphi|^2 \sigma^{-2 \beta}  1_{\sigma > \epsilon}
\leq 2 \beta \int_\Omega \varphi \sigma^{-2 \beta -1} \nabla \varphi \cdot \nabla \sigma 1_{\sigma > \epsilon}
+ \mu^2 
\int_\Omega \varphi^2 \sigma^{-2 \beta}  .
\eea
Noting that $|\nabla \sigma|=1$, applying Cauchy-Schwarz, and rearranging terms, 
we obtain
\bean \label{reuse}
\f12 \int_\Omega |\nabla \varphi|^2 \sigma^{-2 \beta}  1_{\sigma > \epsilon}
\leq 2 \beta^2 \int_\Omega \varphi^2 \sigma^{-2 \beta -2} 
+ \mu^2 
\int_\Omega \varphi^2 \sigma^{-2 \beta} .
\eean
Now using lemma \ref{odd u prop}, as $-2 \beta -2 > -4$,  letting $\epsilon \ri 0$,
the result is proved. $\square$ \\

Interestingly, now that we know that $\sigma^{-\beta} \nabla \varphi$ is in
$L^2(\Omega)$ for $\beta <1$, the estimate in the proof
of lemma \ref{odd u prop} can be improved leading to stronger
decay estimates near the cusp.

\begin{prop}\label{stronger est}
Let $\varphi$ be as in lemma \ref{odd u}.
$\sigma^{- \alpha} \varphi$ and $\sigma^{1- \alpha} \nabla \varphi$ 
are in $L^2(\Omega)$ for any $\alpha <3$.
\end{prop}
\textbf{Proof:}
The proof of lemma \ref{odd u prop} indicates that 
 $\int_{\Omega_R} |\varphi|^2 \leq  16  R^4 \int_{\Omega_R} 
|\p_{x_2} \varphi|^2 $. Thus for $\beta <1$,
\bea
\int_{\Omega_R} |\varphi|^2 \leq 16 R^4 
\int_{\Omega_R} \sigma^{2 \beta}
|\sigma^{-\beta}\p_{x_2} \varphi|^2
\leq 32 R^{4+2\beta}  \int_{\Omega_R} 
|\sigma^{-\beta}\p_{x_2} \varphi|^2,
\eea
for all $R$ small enough.
Following the calculation in the proof of lemma \ref{odd u prop}, we now have,
\bea 
\int_{\Omega_1} |\varphi(x_1, x_2)|^2 ((x_1+1)^2 + x^2)^{- \alpha} 
dx_1 dx_2 \leq
 C \sum_{n=0}^\infty \f{2^{(2n+2) \alpha}}{(2^{n})^{4+2 \beta}},
\eea
which is finite if $\alpha < 2 + \beta$, thus
$\sigma^{-\alpha} \varphi$ is in $L^2(\Omega)$ if
$\alpha<3$.
In other words, $\sigma^{- \alpha} \varphi$ is in $L^2(\Omega)$ if
 $\alpha <3$.
We now rewrite estimate \eqref{reuse} with a parameter $\gamma$,
\bea
\f12 \int_\Omega |\nabla \varphi|^2 \sigma^{-2 \gamma}  1_{\sigma > \epsilon}
\leq 2 \gamma^2 \int_\Omega \varphi^2 \sigma^{-2 \gamma -2} 
+ \mu^2 
\int_\Omega \varphi^2 \sigma^{-2 \gamma} .
\eea
For $\gamma < 2$, the right hand side is finite and we let $\epsilon $
tend to zero. 
$\square$ \\

\begin{lem} \label{ext u prop}
Let $\varphi$ be as in lemma \ref{odd u}.
$\varphi$ can be extended to a function in $H^1(\RR^2)$ which is odd
in $x_2$.
\end{lem}
\textbf{Proof:}
Let $\eta$ be a function in $C^\infty_c(\RR^2)$ equal to 1 in 
a neighborhood of $(-1,0) $ and zero if $ x_1>0$.
Set $\varphi' = \varphi \eta$.
Let 
\bea H^1_{\sigma^{\f12}} (\RR^2):= \{  
\sigma^{\f12}  v, \sigma^{\f12}\nabla v \in L^2(\RR^2)\},
\eea
endowed with its natural norm,
\bea
\|\sigma^{\f12}  v\|_{L^2(\RR^2)} +
\|\sigma^{\f12} \nabla v \|_{L^2(\RR^2)}.
\eea
Given that the  cusp of $\Omega$ at $(-1,0)$ is  of order 2,
there is a continuous extension operator from $H^1(\Omega)$
to $H^1_{\sigma^{\f12}} (\RR^2)$: 
see \cite{maz2013sobolev}, section 1.5.4.
Proposition \ref{stronger est} implies that $\sigma^{-\f32} \varphi'$
is in $H^1(\Omega)$.  Let $E$ be its extension to 
$H^1_{\sigma^{\f12}}(\RR^2) $. We may assume that the support of
$E$ is compact. 
By definition of  $H^1_{\sigma^{\f12}} (\RR^2)$, 
$\sigma^{\f12 } E$ and $\sigma^{\f12 } \nabla E$ are in $L^2(\RR^2)$.
Thus $\sigma^{\f32 }  E $ is in $L^2(\RR^2)$.
But  $\nabla (\sigma^{\f32}E) = \sigma^{\f32} \nabla E + \f32
\sigma^{\f12} E \nabla \sigma $ is also in $L^2(\RR^2)$, as 
$\nabla \sigma$ is bounded.
Thus $\sigma^{\f32}E$ is an extension of $\varphi '$ which is in 
$H^1(\RR^2)$. The cusp at $(1,0)$ can be handled likewise:
we extend $\varphi (1 -\eta)$ to a function in $H^1(\RR^2)$.
Adding the two extensions, we obtain a function $\tilde{E}$
in $H^1(\RR^2)$ equal to $\varphi$ in $\Omega$.
Finally, $\f12 (\tilde{E}(x_1, x_2)) - \tilde{E}(x_1, -x_2))$
is in $H^1(\RR^2)$, is odd in $x_2$, and equals 
$\varphi$ in $\Omega$, since $\varphi(x_1,-x_2) = -\varphi(x_1,x_2)$.
$\square$ \\\\

\textsl{Construction of the related counterexample for the unique solvability of the crack inverse problem.}\\
Let $\varphi$ be as in lemma \ref{odd u}.
According to lemma  \ref{odd u} 
the function $\varphi^+$ defined by 
$\varphi^+(x_1, x_2) = \varphi(x_1, x_2)$
if $x_2 >0, (x_1,x_2) \in \Omega$, $\varphi^+(x_1, x_2) =0$
if $x_2 <0, (x_1,x_2) \in \Omega$, is still in $H^1(\Omega)$
and according to lemma \ref{ext u prop}, $\varphi^+$ can be extended 
to a function in $H^1(\RR^2)$ which is zero if $x_2<0$.
This
explains why the trace of $\varphi$ restricted  to
$\p \Omega^+ = \{ (x_1,x_2) \in \p \Omega: x_2>0 \}$ is
in $\tilde{H}^{\f12} (\p \Omega^+)$.\\
Define $\Gamma_1 = \{ (x_1, x_2): (x_1,x_2+2) \in \p \Omega^+ \}$,
$\Gamma_2 = \{ (x_1, x_2): (x_1,x_2+2) \in \p \Omega, x_2+2<0\}$,
and
$g_i(x_1,x_2) = \varphi(x_1,x_2+2), (x_1,x_2) \in \Gamma_i$, $i=1,2$.
We now know that
$g_i$ is  in  $\tilde{H}^{\f12} (\Gamma_i)$.
Let 
 $u^i$ be the unique solution to  the two dimensional analog
of system (\ref{BVP1}-\ref{Decay1}) with $\Gamma_i$ in place of $\Gamma$, where 
the continuous normal vector on $\Gamma_i$ is such that $n \cdot e_2 >0$,
$g^i$ in  $\tilde{H}^{1/2}(\Gamma_i)$ is in place of $g$,
 $k^2 = \mu^2$ as in lemma \ref{odd u},
and the solutions are sought in the functional space
\bea
{\cal V}_2 := \{ v \in H^1_{loc}(\RR^{2-}\setminus \ov{\Gamma}): \f{v}{\sqrt{1+r^2}}, 
\f{\nabla v}{\sqrt{1+r^2}}, \f{\p v}{\p r} - i \mu v \in L^2(\RR^{2-}\setminus \ov{\Gamma})
\}.
\eea
Let $u= u^1 - u^2$.
We can show that $u$ is zero outside $\Omega - (0,0,2)$ by repeating the same argument
as in the end of section \ref{counter-example}.
Note that $g_i$ has full support in $\Gamma_i$: since 
$\f{\p \varphi}{\p n}=0$ on $\p \Omega$,
this is again due to the corollary of  theorem 1 in \cite{barcelo1988weighted}.\\
Finally, by using the rescaling $(x_1, x_2) \ri (s x_1 , s x_2)$, we can find a similar
counterexample 
for any $k^2>0$.
Replacing $f$ by 
  $f_a(x_1)=a(x_1-1)^2(x_1+1)^2 - 2$, $x_2 
\in (-1,1)$, where $0<a<1$ instead of $f$, 
the same argument as above can be repeated using 
the domain with cusps $\{ -f_a(x_1)< x_2 < f_a(x_1): -1< x_2 <1\}$.
Picking $a$ small can make the corresponding curves $\Gamma_1$
and $\Gamma_2$ arbitrarily flat.
Note that $\Gamma_1$ and $\Gamma_2$ can be extended 
by line segments while preserving their $C^1$ regularity. 

\subsection{A counterexample  in case of cracks
that are nearly flat, the three-dimensional case}\label{counter-example3}
The previous two-dimensional counterexample can serve as 
an inspiration for 
constructing a three-dimensional counterexample with cylindrical
symmetries. 
We have to implement relevant modifications to the two dimensional
functional spaces that we are using.
Let $(r,\theta, x_3)$ denote the cylindrical coordinates in $\RR^3$.
Since we will focus on functions that are
independent of $\theta$, for any open subset
$U$ of $\RR^2$  where the points in $U$ are denoted by $(r,x_3)$,
define the space $L^2_r(U)$ of measurable functions $u$ in $U$ such that 
$\int_U u^2 r drdx_3 < \infty$ equipped with the norm
\bea
\| u \|_{L^2_r(U)} = (\int_U u^2 r dr  dx_3)^\f12 .
\eea
Similarly, let $H^1_r(U)$ be the subspace of functions $u$ in
${L^2_r(U)} $ that have weak derivatives in  ${L^2_r(U)} $.
$H^1_r(U)$ will be equipped with the norm 
\bea
\| u \|_{H^1_r(U)} = \| u \|_{L^2_r(U)} +
(\int_U |\nabla u|^2 r dr  dx_3)^\f12 ,
\eea
where $\nabla u =(\p_r u, \p_{x_3} u)$.
Next, we define the open subset of $\RR^2$,
\bea
\Omega_0 :=\{ (r,x_3):  0<r<1, -f(r) < x_3 < f(r) \},
 \eea
where $f$ is the same  as in \eqref{omega def}, and the open
subset of $\RR^3$,
\bea
\Omega_{\RR^3} :=\{ (r,\theta, x_3):  0<r<1, 0\leq \theta <2 \pi,  -f(r) < x_3 < f(r) \}.
 \eea
The shape in figure \ref{omegaRdomain} is then a cross-section 
of $\Omega_{\RR^3}$ by any plane containing the $x_3$ axis.

\begin{lem}\label{comp embed3}
$H^1_r(\Omega_0) $ is compactly embedded in $L^2_r(\Omega_0)$.
\end{lem}
\textbf{Proof:}
We set $U_1 = \{ (r,x_3) \in \Omega_0:  r < \f34 \}$, 
$U_2 = \{ (r,x_3) \in \Omega_0 : r > \f14 \}$.
Rotating $U_1$ about the $x_3$ axis, we obtain a bounded Lipschitz
domain $U_{1,\RR^3}$ in $\RR^3$. As $H^1(U_{1,\RR^3})$
is compactly embedded in $L^2(U_{1,\RR^3})$, it follows 
that $H^1_r(U_1) $ is compactly embedded in $L^2_r(U_1)$.
In $U_2$,  $\f14 < r <1$, thus the spaces $L^2_r(U_2)$
and $L^2(U_2)$ are the same, as  well as  the spaces 
 $H^1_r(U_2)$
and $H^1(U_2)$. It now follows from lemma \ref{comp embed} that
$H^1_r(U_2) $ is compactly embedded in $L^2_r(U_2)$.\\
Finally, fix a smooth function $p:(0, \infty) \ri \RR$ such that $0 \leq p \leq 1$,
$p(r) =1$, if $  r < \f14$, $p(r) =0$, if $  r > \f34$.
For $u$ in $H^1_r(\Omega_0) $, write
$u = p u + (1-p) u$ and the result follows. 
$\square$\\

\begin{lem} \label{odd u3}
There exists a  function $\varphi$ in $H^1(\Omega_{\RR^3})$ and $\mu>0$ such that 
$\int_\Omega \varphi^2=1$,
$(\Delta+\mu^2) \varphi =0$, $\f{\p \varphi}{\p n} =0$, almost everywhere
on $\p \Omega_{\RR^3}$,
$\varphi$ does not depend on $\theta$,
 and $\varphi(r, -x_3) = - \varphi(r, x_3)$, for $(r,x_3)$
  in $\Omega$.
\end{lem}
\textbf{Proof:}
We introduce the two subspaces of $H^1_r(\Omega_0)$
\bea
H^{1,\pm}_r(\Omega_0):= \{v \in H^1_r(\Omega_0): v(r,-x_3) = \pm 
v(r,x_3), (r,x_3) 
\in \Omega_0 \}.
\eea
They are orthogonal to one another for the natural inner product
in $H^1_r(\Omega_0)$. 
Let
\bea
\mu^2= \inf_{v \in H^{1,-}_r(\Omega_0), \| v \|_{L^2_r(\Omega_0)} =1} \int_{\Omega_0} |\nabla v|^2 rdrdx_3.
\eea
Again, 
a minimizing sequence for this inf converges to some $\varphi_0$ in
$H^{1,-}_r(\Omega_0)$ such that for all $\phi$ in $\in H^{1,-}_r(\Omega_0)$
\bea
\int_{\Omega_0} ( \nabla \varphi_0 \cdot \nabla \phi - 
\mu^2 \varphi_0 \phi)  r dr dx_3 =0,
\eea
 and by orthogonality, this holds for
all $\phi$ in $\in H^{1}_r(\Omega_0)$.
We then set the extension $\varphi(r,\theta,x_3) =
\varphi_0(r,x_3) $ to obtain a function in $H^1(\Omega_{\RR^3})$ 
which has the desired properties.
$\square$\\

\begin{prop} \label{stronger est3}
Let $\varphi_0$ be as in the proof of lemma 
\ref{odd u3}.
Let $\sigma$ be the weight $\sigma(r,x_3) = \sqrt{(r-1)^2 + x_3^2 }$.
Then $\sigma^{- \alpha} \varphi_0$ and $\sigma^{1- \alpha} \nabla \varphi_0$ 
are in $L^2_r(\Omega_0)$ for any $\alpha <3$.
\end{prop} 
\textbf{Proof:}
Just as its analog in  the purely two dimensional case of section
\ref{counter-example2}, this proposition is proved in two steps.
In the first step we prove it for $\alpha < 2$. 
Let $R \in (0,1)$ and 
$$\Omega_{0,R} :=\{  (r,x_3) \in \Omega_0: 1-R < r < 1 \}. $$
Since in $\Omega_{0,R}$, $r$ is bounded above by 1, 
the proof is nearly identical to that of lemmas \ref{odd u prop}
and \ref{grad u prop}. Next, the decay estimates can  be improved since
\bea
\int_{\Omega_{0,R}} |\varphi_0|^2 \leq 16 R^4 
\int_{\Omega_{0,R}} \sigma^{2 \beta}
|\sigma^{-\beta}\p_{x_3} \varphi_0|^2 rdrdx_3
\leq 32 R^{4+2\beta}  \int_{\Omega_{0,R}} 
|\sigma^{-\beta}\p_{x_3} \varphi_0|^2 rdrdx_3,
\eea
for all $R$ small enough,
and the rest of the proof follows as in the proof of proposition 
\ref{stronger est},
with minor adjustments due to the different surface are element in the present 
case.
$\square$\\

\begin{lem} \label{ext u prop3}
Let $\varphi$ be as in lemma \ref{odd u3}.
$\varphi$ can be extended to a function in $H^1(\RR^3)$ which is
 odd in $x_3$.
\end{lem}
\textbf{Proof:}
We use again the  smooth function $p:(0, \infty) \ri \RR$ such that $0 \leq p \leq 1$,
$p(r) =1$, if $  r < \f14$, $p(r) =0$, if $  r > \f34$.
It suffices to  find separate extensions of $p \varphi$ and 
$(1-p) \varphi$ from $H^1(\Omega_{\RR^3})$ to $H^1(\RR^3)$.\\
It is clear that $p \varphi$ can be extended from 
$H^1(\Omega_{\RR^3})$ to $H^1(\RR^3)$
since the bounded domain 
$\Omega_{\RR^3} \cap \{|r|<\f34 \} $ is Lipschitz regular.\\
The function $\sigma^{- \f32} (1-p) \varphi_0$ is in 
$H^1_r(\Omega_0)$ thanks to proposition \ref{stronger est3}.
As it is zero if $r < \f14$, it is also in $H^1(\Omega_0)$,
so using the same argument as in lemma \ref{ext u prop},
it can be extended to a function $E_0$ in $H^1_{\sigma^{\f12}}(\RR^2)$.
We can assume $E_0$ to be compactly supported. 
Then, just as in the proof of lemma \ref{ext u prop}, we claim
that $\sigma^{\f32} E_0$ is in $H^1(\RR^2)$. It is also
in $H_r^1(\RR^2)$, since it has compact support.
Let $\psi$ be function in $H^1(\RR^3)$ which is independent of the polar
angle $\theta$ and whose cross-section is  $\sigma^{\f32} E_0$.
Note that $\psi$ equals $(1-p) \varphi$ in $\Omega_{\RR^3}$.
Finally, we add the extensions of $p \varphi$  and  $(1-p) \varphi$.
We can ensure that this extension is odd in $x_3$ by only retaining its odd
part, as in the end of the proof of lemma \ref{ext u prop}.
$\square$\\\\

\textsl{Construction of the related counterexample to the unique solvability of the crack inverse problem - the three-dimensional case.}\\
We pretty much follow the construction in the two-dimensional 
case demonstrated at the end of section \ref{counter-example2}.
For $\varphi$ be as in lemma \ref{odd u3}, set
$\varphi^+= \varphi$
if $x_3 >0, x \in \Omega_{\RR^3}$, $\varphi^+ =0$
if $x_3 <0, x \in \Omega_{\RR^3}$.
$\varphi^+$ is still in $H^1(\Omega_{\RR^3})$
according to lemma  \ref{odd u3} 
and according to lemma \ref{ext u prop3}, $\varphi^+$ can be extended 
to a function in $H^1(\RR^3)$ which is zero if $x_3<0$.
This
explains why the trace of $\varphi$ restricted  to
$\p \Omega_{\RR^3}^+ = \{ x \in \p \Omega_{\RR^3}: x_3> 0\}$ is
in $\tilde{H}^{\f12} (\p \Omega_{\RR^3}^+)$.\\
Define $\Gamma_1 = \{ (r, \theta , x_3): (r, \theta ,x_3+2)
 \in \p \Omega_{\RR^3}^+ \}$,
$\Gamma_2 = \{ (r, \theta , x_3): (r, \theta ,x_3+2) \in
 \p \Omega_{\RR^3}, x_3+2<0\}$,
and
$g_i(r, \theta , x_3) = \varphi(r, \theta , x_3+2), (r, \theta , x_3) \in \Gamma_i$, $i=1,2$.
We now know that
$g_i$ is  in  $\tilde{H}^{\f12} (\Gamma_i)$.
Let 
 $u^i$ be the unique solution to  
system (\ref{BVP1}-\ref{Decay1}) with $\Gamma_i$ in place of $\Gamma$, where 
the continuous normal vector on $\Gamma_i$ is such that
 $n \cdot e_3>0$,
$g^i$ in  $\tilde{H}^{1/2}(\Gamma_i)$ is in place of $g$,
 $k^2 = \mu^2$ as in lemma \ref{odd u3}.
Let $u= u^1 - u^2$.
We can show that $u$ is zero outside $\Omega - (0,0,2)$ by repeating the same argument
as in the end of section \ref{counter-example}.
Note that $g_i$ has full support in $\Gamma_i$: since $\f{\p \phi}{\p n}=0$ on $\p \Omega$,
this is again due to the corollary of  theorem 1 in \cite{barcelo1988weighted}.\\

\section{The crack inverse problem for pressure waves in half-space: uniqueness of solutions except for a
discrete set of frequencies} \label{unique discrete}
We have covered examples where  two cracks $\Gamma_1$ and
$\Gamma_2$ are such that $\RR^3\sm \ov{\Gamma_1 \cup \Gamma_2}$
has more than one connected component and there are some $k$ and $g_1$ and $g_2$
such that if 
$u^1$ and $u^2$ solve problem (\ref{BVP1}-\ref{Decay1}) with respective 
forcing terms $g_1$ and $g_2 $,  we have  that $u_1 = u_2  $ on $\{x: x_3 = 0\}$
even if $\Gamma_1 \neq \Gamma_2$ and $g_1, g_2 $ have full support in
$\Gamma_1, \Gamma_2$.
By contrast, we now prove that under a general requirement on the geometry
of $\Gamma$
this can only occur for a discrete set of
frequencies, in other words,  uniqueness for the inverse problem holds for wavenumbers
$tk$, except possibly for a discrete set of scaling parameters $t$.
 Let ${\cal P}$ be the class of Lipschitz open surfaces that are finite unions of polygons. 
Since any polygon is a finite union of triangles, ${\cal P}$ is  equivalently the class of 
Lipschitz open surfaces that can be obtained from a finite triangulation. 
Fix $k$ in $L^\infty(\RR^{3-})$ satisfying conditions $(H1)$ through $(H3)$. 
For $t>0$ consider the problem
\bean
        (\Delta + t^2k^2 )u=0\text{ in }\doubleR^{3-}\sm \ov{\Gamma},  \label{BVP12}     \\
        \p_{x_3} u=0\text{ on  the surface } x_3=0,  \label{BVP22}     \\
     \aaa \frac{\p u}{\p n} \bbb = 0 \mbox{ across }\Gamma,   \label{BVP32}\\
      \aaa u \bbb
			=g \mbox{ across }\Gamma,  \label{BVP42}\\
    u \in {\cal V}_t \label{Decay12},
\eean
where
\bean
{\cal V}_t := \{ v \in H^1_{loc}(\RR^{3-}\setminus \ov{\Gamma}): \f{v}{\sqrt{1+r^2}}, 
\f{\nabla v}{\sqrt{1+r^2}}, \f{\p v}{\p r} - i t k_0 v \in L^2(\RR^{3-}\setminus \ov{\Gamma})
\}. 
\eean

\begin{thm}
\label{InverseProblemResult2}
For $i=1,2$,
    let $\Gamma_i$   be a Lipschitz open surface
			of  class ${\cal P}$
			such that its closure is
			in
		 $\doubleR^{3-}$, let 
 $u^i$ be the unique solution to  (\ref{BVP12}-\ref{Decay12}) with $\Gamma_i$ in place of $\Gamma$ and the jumps $g^i$ in  $\tilde{H}^{1/2}(\Gamma_i)$ in place of $g$.  
Assume that $g_i$ has full support in $\Gamma_i$.
	There is a discrete set ${\cal D}$ in $[0, \infty)$
	such that if $t \notin {\cal D} $, 
	if  $V$ is a non-empty relatively open subset of $\{x: x_3=0\}$, 
	and 
	if $u^1=u^2$ on $V$, then
	$\ov{\Gamma_1}=\ov{\Gamma_2}$ and $g^1=g^2$ almost everywhere.	
\end{thm}

\textbf{Remark}:
If $\Gamma_1$ and $\Gamma_2$ are allowed to be general Lipschitz domains,
then a connected component $\Omega_i$ of $\RR^{3-} \setminus \ov{\Gamma_1 \cup \Gamma_2}$ 
may be one of those so called bad domains discussed in chapter 2 section 5 of 
\cite{maz1998differentiable} (cusps may not be of power type). In particular, 
$H^1(\Omega_i)$ may not be compactly embedded in $L^2(\Omega_i)$ which significantly complicates the study of Neumann eigenvalues. \\
\textbf{Proof} of theorem \ref{InverseProblemResult2}:
Since $\Gamma_1, \Gamma_2$ are  Lipschitz open surfaces
		of  class ${\cal P}$,
$\RR^{3-} \setminus \ov{\Gamma_1 \cup \Gamma_2}$ has a finite
	number of connected components $\Omega_i$, $i=1,..,N$.
	If $N=1$, the result holds due to theorem \ref{InverseProblemResult}.
	If $N \geq 2$, we may assume that $\Omega_1$ is unbounded, while $\Omega_i$,
	$i \geq 2$ is bounded. By construction, each $\Omega_i$, $i \geq 2$, is a bounded Lipschitz domain.
	Thanks to   conditions $(H1-H2)$, we can define on $H^1(\Omega_i)$, $i \geq 2$,
	the  inner product,
	\bea
	{\cal B}(v,w) = \int_{\Omega_i} \nabla v \cdot \nabla w + k^2 v w,
	\eea
	which is equivalent to the the natural inner product.
	Similarly, we  define on $L^2(\Omega_i)$,
	\bea
	b(v,w) = \int_{\Omega_i}  k^2 v w .
	\eea
	By Lax-Milgram's theorem, for any continuous linear functional $F$ on $H^1(\Omega_i)$,
	there is a unique $TF$ in $H^1(\Omega_i)$ such that 
	${\cal B}(TF,w) = Fw$, for all $w$ in $H^1(\Omega_i)$, and $T$ is a continuous linear operator.
	Since $\Omega_i$ is a Lipschitz domain,  $H^1(\Omega_i)$  is compactly embedded
	in $L^2(\Omega_i)$, thus we can define a compact operator $K$ from $H^1(\Omega_i)$
	to its dual   mapping $\phi$ to the functional $w \ri b(\phi, w)$.
	$TK$ is now a compact operator from $H^1(\Omega_i)$ to $H^1(\Omega_i)$ and satisfies
	${\cal B}(TKv,w)=b(v,w)$, for all $v, w$ in $H^1(\Omega_i)$.
	It is clear that $TK$ is symmetric, positive, and definite, thus its eigenvalues form
	a decreasing sequence of positive numbers $\alpha_n$ converging to zero.
	If $\alpha$ is not
	an eigenvalue, then $TK - \alpha I$ is invertible.
	Let ${\cal N}_i$ be the set of these eigenvalues.
	If $\alpha \notin {\cal N}_i$ and $v$ is in $H^1(\Omega_i)$, 
	if for all $w$ in $H^1(\Omega_i)$,  ${\cal B}((TK- \alpha I)  v, w) =0$, then $v=0$.
	As ${\cal B}((TK- \alpha I)  v, w) = b(v,w) - \alpha {\cal B}(v,w)$,
	it follows that if $(\Delta + (\f{1}{\alpha} - 1)k^2) v= 0$ in $\Omega_i$ and $\f{\p v}{\p n} = 0$
	on $\p \Omega_i$, then $v =0$.
	We then set 
	\bea
	 {\cal D}_i = \{  \sqrt{\f{1}{\alpha_n} - 1}: n \geq 1  \} ,
	\eea
	and ${\cal D}= \cup_{i=2}^N {\cal D}_i$.
	Assume that $t \notin {\cal D}$.
	Let $U=\bigcup_{i=1}^N \Omega_i$ and set $u= u^1 - u^2$ in $U$.
Since $(\Delta + t^2 k^2) u =0$ in $U$ and $ u = \p_{x_3} u = 0$ on $V$, $u$ can be extended by zero
 to an open set $U'$ of $\RR^3 $ such that  $\Omega_1 \subset U'$, $U'$ is connected, 
$U' \cap \{x: x_3>0\}$ is non-empty, $u$ is in $H^2_{loc} (U')$ and satisfies $(\Delta +t^2 k^2) u  =0$ in 
$U'$, where we can set $k=0$ in  $U' \cap \{x: x_3>0\}$. 
 By the unique continuation property (corollary of  theorem 1 in \cite{barcelo1988weighted}),
 $u$ is zero in $U'$. Next, using \eqref{BVP32},  we find that 
$\frac{\p u}{\p n}=0$ almost everywhere on $\p \Omega_i$ for $i \geq 2$: as $t \notin {\cal D}$,
it follows that $u$ is also zero in $\Omega_i$.
To finish the proof, we just need to carry out the same argument as in the proof 
of theorem \ref{InverseProblemResult}. 
$\square$.

\end{document}